\documentclass[a4paper,12pt]{amsart}

\usepackage[latin1]{inputenc}
\usepackage{amsmath}
\usepackage{amsthm}
\usepackage{amssymb}
\usepackage{graphicx}
\usepackage{bbm}
\usepackage{microtype}
\usepackage{xcolor}
\usepackage{hyperref}
\usepackage{cleveref}

\newcommand*\diff{\mathop{}\!\mathrm{d}}
\newcommand{\suchthat}{\;\ifnum\currentgrouptype=16 \middle\fi|\;}
\newcommand{\R}{\mathbb{R}}
\renewcommand{\S}{\mathcal{S}}
\newcommand{\h}{\widehat}
\renewcommand{\l}{\left(}                                  
\renewcommand{\r}{\right)}                               %
\renewcommand{\epsilon}{\varepsilon}
\renewcommand{\phi}{\varphi}
\newcommand{\lv}{\left \Vert}
\newcommand{\rv}{\right \Vert}

\theoremstyle{plain}
\newtheorem{thm}{Theorem}
\newtheorem{cor}[thm]{Corollary}

\newtheorem{lemma}[thm]{Lemma}

\theoremstyle{definition}
\newtheorem{defn}[thm]{Definition}
\newtheorem{rem}[thm]{Remark}

\title[A logarithmic uncertainty principle for functions with radial symmetry]{A logarithmic uncertainty principle for functions with radial symmetry}
\author{Jacopo Bellazzini}
   \address{Jacopo Bellazzini \newline Dipartimento di Matematica  \\ Universit\`a di Pisa \\t Largo B. Pontecorvo 5, 56100 Pisa, Italy}%
\author{Matteo Nesi}
\address{Matteo Nesi \newline Scuola Normale Superiore \\ Piazza dei Cavalieri 7, 56100 Pisa, Italy}

\begin{document}

\thanks{J. B. is partially supported by project PRIN 2020XB3EFL by the Italian Ministry of Universities and Research and by  the University of Pisa, Project PRA 2022 11.}

\begin{abstract}
In this note, we prove a new uncertainty principle for functions with radial symmetry by differentiating a radial version of the Stein-Weiss inequality.
The difficulty is to prove the differentiability in the limit of the best constant that, unlike the general case, it is not known.
\end{abstract}

\maketitle

\section{Introduction}
In this note, we will be differentiating inequalities, so we start by summarizing the general reasoning scheme. Let us suppose that we are given a family of inequalities 
\[
	F_t \leq k_t G_t
\]
parametrized by $t>0$. We choose to keep separate the 
families $F_t, G_t$ from the family $k_t$, since in later applications 
$F_t, G_t$ will be given by some functionals computed on a fixed function $f$, while $k_t$ are the (optimal) constants in the inequalities.
Suppose further that
\begin{itemize}
\item the functions $F_t,G_t,k_t$ are actually defined and differentiable for $t\geq0$ ;
\item $F_0=k_0G_0$.
\end{itemize}
Then we can write
\[
	\frac{F_t-F_0}{t} \leq \frac{k_tG_t-k_0G_0}{t}
\]
and passing to the limit for $t\to0^+$ we get
\begin{equation}\label{eq:disegprinc}
	\partial_t F_t |_{t=0} \leq G_0 \partial_t k_t |_{t=0}+k_0 \partial_t G_t |_{t=0} .
\end{equation}

These computations will be carried out, in each case, for functions in a suitable class which allows differentiation under the integral sign; in our case it is convenient to consider the Schwartz class $\S(\R^d)$, which is made of regular and decaying functions and is well-behaved with respect to the Fourier transform.

For a first example, we start with the Hardy-Littlewood-Sobolev inequality
\[
	\left| \int_{\R^d} \int_{\R^d} \frac{f(x)g(y)}{|x-y|^\lambda} \diff x \diff y \right| \lesssim ||f||_p ||g||_q,
\]
which, in the case $f=g\geq0$ ($p=q$ and $0<\lambda<d$), becomes
\[
	 \int_{\R^d} \int_{\R^d} \frac{f(x)f(y)}{|x-y|^\lambda} \diff x \diff y  \leq k_{\lambda} ||f||^2_{\frac{2d}{2d-\lambda}},
\]
where
\[
	k_{\lambda}=\pi^{\frac{\lambda}{2}} \frac{\Gamma(\frac{d}{2}-\frac{\lambda}{2})}{\Gamma(d-\frac{\lambda}{2})}\left(\frac{\Gamma(\frac{d}{2})}				{\Gamma(d)}\right)^{-1+\frac{\lambda}{d}}
\]
 and the inequality  becomes an \emph{equality} when $\lambda=0.$ By observing that 
\[
	\frac{\diff}{\diff \lambda} |x-y|^{-\lambda} |_{\lambda=0} = -\log|x-y|
\]
and following the above reasoning scheme, we get the logarithmic Hardy-Littlewood-Sobolev inequality (cf. \cite{CarlenLoss}):
\[
	\int_{\R^d} \int_{\R^d} \! -\log|x-y| \ f(x) f(y) \diff x \diff y \leq \frac{1}{d} \int_{\R^d} f \log f \diff x + C_0,
\]
for regular $f\geq0$ with $||f||_1=1$, where $C_0 = \partial_{\lambda} k_{\lambda} |_{\lambda=0}$.

Also, a quantity similar to the LHS of the previous inequality can be estimated differentiating the sharp Young inequality: as in \cite{LiebLoss} (ch. 8, p. 223), we find, for any $a>0$, the logarithmic Sobolev inequality 
in $\R^d$
\[
	\int_{\R^d} |f|^2 \log |f| \diff x \leq \frac{a^2}{\pi} \int_{\R^d} |\nabla f|^2 \diff x -d(1+\log a),
\]
which holds for every $f\in H^1(\R^d)$ with $\lv f \rv_2=1$.

Concerning the Fourier transform, a classical inequality is the \\
Hausdorff-Young inequality for $1\leq p \leq 2$
\[
	||\h{f}||_{p'} \lesssim ||f||_p,
\]
from which Beckner \cite{Beckner:Fourier} derived a logarithmic uncertainty principle for $f\in \S(\R^d)$ with $\lv f \rv_2=1$
\[
	\int_{\R^d} |f|^2 \log |f| \diff x + \int_{\R^d} |\h{f}|^2 \log |\h{f}| \diff \xi
			\leq \frac{d}{2} (\log(2)-1)<0.
\]
So far, we have examples in which we get inequalities involving logarithmic terms, in the form of $\log(f)$ or the kernel $\log|x-y|$. Now we recall three families of inequalities which give also logarithmic weights in physical $(\log|x|)$ and Fourier $(\log|\xi|)$ domains. These are 
Sobolev, Hardy and Stein-Weiss inequalities, the latter being a generalization of the first two "endpoint" inequalities. The results of the differentiation are summarized as follows: the Sobolev and Hardy inequalities give, respectively, ($\lv f \rv_2=1$)  \cite{Beckner:Pitt}
\begin{align*}
 ||f||_{\frac{2d}{d-2s}}&\lesssim ||f||_{\dot{H}^s} \ \longrightarrow \ \frac{2}{d}\int_{\R^d} |f|^2 \log |f| \diff x 
		\leq \int_{\R^d} |\hat{f}|^2\log |\xi| \diff \xi  + c \\
\lv |\cdot|^{-s} f \rv_2 &\lesssim ||f||_{\dot{H}^s} \ \longrightarrow 
	\quad -\int_{\R^d} |f|^2 \log|x| \diff x  \leq \int_{\R^d} |\hat{f}|^2 \log |\xi| \diff \xi + c
\end{align*}
and, by interpolating, the Stein-Weiss inequality
\[
	 \lv |\cdot|^{-\beta} f \rv_p \lesssim ||f||_{\dot{H}^s} 
\]
yields (cf. \cref{dSW}) 
\begin{equation}\label{idSW}
\begin{split}  
	\frac{2(1-t)}{d}\int_{\R^d}& |f|^2 \log |f| \diff x  - t \int_{\R^d} |f|^2 \log|x| \diff x \leq \\
		&\leq \int_{\R^d} |\h{f}|^2 \log |\xi| \diff \xi  + c 
\end{split}
\end{equation}
for all $t\in[0,1]$. In the above inequalities, $c$ is a constant that differs from line to line. 

\medskip

The aim of the note is to prove an improvement of the logarithmic Stein-Weiss inequality (\cref{idSW}) when restricted to radially symmetric functions. We denote by $\dot{H}_{rad}^s(\R^d)$ the space of radial functions which belong to $\dot{H}^s(\R^d)$. The starting point is an improvement of Stein-Weiss inequality for radial functions, which is due to Rubin.
\begin{thm}[Rubin inequality, \cite{Rubin}, \cite{BFV}, \cite{Bellazzini}]  \label{Rubin}
Given $f\in \dot{H}_{rad}^s(\R^d)$ with $ 0\leq s < d/2$ and $p,\beta$ with the relations $(1/0 := +\infty)$
\begin{gather*}
2 \leq p \leq \frac{2}{(1-2s)_+} \\
\beta = s+\frac{d}{p}-\frac{d}{2} ,
\end{gather*}
one has the inequality
\[
	\lv \  |\cdot|^{-\beta} f \rv_p \leq C(p,s) ||f||_{\dot{H}^s} .
\] 
\end{thm}

\begin{rem}
The latter theorem is an improvement of Stein-Weiss inequality on radial functions, meaning that the range of parameters $p$ extends further on the right and for big values of $p$ even positive exponents $-\beta$ are allowed.
\end{rem}

We point out that in \cref{Rubin} the optimal constants $C(p,s)$ are not known, but via an interpolation argument we get an estimate which is sufficient to carry out the differentiation, so as to get:

\begin{thm}\label{dRubin}
There exists a constant $c\in\R$ such that, given $f\in \S(\R^d)$ radial and with $\lv f \rv_2=1$
\begin{equation}\label{eq:main}
\begin{split}
	2\int_{\R^d} |f|^2 &\log |f| \diff x   + (d-1) \int_{\R^d} |f|^2 \log|x| \diff x \leq \\
		&\leq \int_{\R^d} |\h{f}|^2 \log |\xi| \diff \xi  + c 
\end{split}
\end{equation}

\end{thm}

\begin{rem} \label{dSW}
The argument in the proof holds, in particular, in the restricted range of parameters $2\leq p \leq 2d/(d-2s)$ or equivalently $0\leq\beta\leq s$. For this range of exponent, the Rubin inequality holds even for non-radial functions (thanks to Stein-Weiss inequality) and the differentiation gives, for any $0\leq t \leq 1$, again 
\[
	\frac{2(1-t)}{d}\int_{\R^d} |f|^2 \log |f| \diff x  \ - \ t \int_{\R^d} |f|^2 \log|x| \diff x
		\leq \int_{\R^d} |\h{f}|^2 \log |\xi| \diff \xi  + c .
\]
Notice that the relevant difference w.r.t. the radial case is only the sign of the constant in front of the second term.
\end{rem}

\begin{rem}
Stein-Weiss inequality can be improved in the radial case (Rubin inequality), while Sobolev and Hardy cannot.
\end{rem}
\subsection{Notations.} 

$L^p(\R^d)$, with $p\in[1,\infty],$  denote the usual Lebesgue spaces, $\dot H^s(\R^d)$ stands for the usual homogeneous Sobolev space, namely the space of tempered distribution $u$ over $\R^d$, the Fourier transform of which belongs $L^1_{loc}(\R^d)$ and satisfies 
\[
	||u||_{\dot H^s(\R^d)}^2=\int_{\R^d} |\xi|^{2s} |\h{u}|^2d\xi <+\infty.
\]
In this paper, we adopt the following definition for the Fourier Transform:
\[
	\mathcal{F}[f](\xi) := \h{f}(\xi) := (2\pi)^{-d/2} \int_{\R^d} f(x) \ e^{-i \xi \cdot x} \diff x,
\]
so that, defining
\[
	\h{|D|f} := |\xi|\h{f},
\]
one has
\[
	\lv |D| f \rv_2 = \lv \nabla f \rv_2.
\]
For references on the properties of homogeneous and nonhomogeneous Sobolev spaces, we refer to \cite{BCD}. 

\section{Proof of main results}

In order to prove the estimate on the growth of the constant of the Rubin inequality near $s=0$ we need an interpolation result due to Stein, so we start with a preliminary definition

\begin{defn}
Let $(M, \mathcal{M},\mu), (N, \mathcal{N},\nu)$ be two $\sigma$-finite measure spaces and $\{T_z\}_{z\in \overline{\Sigma}}$ a family of operators defined on the space of simple functions on $M$ and taking values in the space of measurable functions on $N$. The family of operators is said to be \emph{admissible} if, for every couple of simple functions $f,g$ defined on $M,N$ respectively, the function
\[
	z\mapsto \langle \overline{T_z f}, g\rangle := \int_N T_z f \ g \diff \nu 
\]
is holomorphic in $\Sigma:= (0,1)\times i\R\subset \mathbb{C}$ and continuous up to the boundary and for some $a<\pi$
\[
	\sup_{z\in \overline{\Sigma}} \ e^{-a |y|} \log |\langle T_z f, g \rangle | <+\infty.
\]
\end{defn}

\begin{thm}[Stein interpolation, \cite{Stein},\cite{Stein:Weiss}] 
Given $p_0, p_1, q_0, q_1 \in [1,+\infty]$ and an admissable family of operators $\{T_z\}_{z\in\overline{\Sigma}}$ satisfying 
\[
	||T_{iy}||_{L^{p_0}\to L^{q_0}}\leq M_0 (y) , \quad ||T_{1+iy}||_{L^{p_1}\to L^{q_1}}\leq M_1 (y)
\]
and the growth condition
\[
	\sup_{y\in\R} \ e^{-a|y|} \log M_j(y) < +\infty \quad \text{ for } j=0,1
\]
with $a<\pi$, it holds
\[
	||T_t||_{L^{p_t} \to L^{q_t}} \leq M_t ,
\]
for $0<t<1$, where
\begin{gather*}
	M_t :=  \int_\R \frac{M_0(y)}{\cosh(\pi y) -\cos(\pi t)} + \frac{M_1(y)}{\cosh(\pi y) +\cos(\pi t)} \diff y , \\
	\frac{1}{p_t} := \frac{1-t}{p_0}+\frac{t}{p_1}, \\
	\frac{1}{q_t} := \frac{1-t}{q_0}+\frac{t}{q_1}.
\end{gather*}
\end{thm}

We are interested in the case where the functions $M_0, M_1$ can be chosen to be bounded; the previous theorem becomes:

\begin{cor}[Stein interpolation - weak version] \label{sw}
Given $p_0, p_1, q_0, q_1 \in [1,+\infty]$ and an admissable family of operators $\{T_z\}_{z\in\overline{\Sigma}}$ satisfying 
\[
	||T_{iy}||_{L^{p_0}\to L^{q_0}}\leq M_0 , \quad ||T_{1+iy}||_{L^{p_1}\to L^{q_1}}\leq M_1 ,
\]
 it holds
\[
	||T_t||_{L^{p_t} \to L^{q_t}} \leq M_t ,
\]
for $0<t<1$, where
\begin{gather*}
	M_t :=  M_0^{1-t} M_1^t , \\
	\frac{1}{p_t} := \frac{1-t}{p_0}+\frac{t}{p_1}, \\
	\frac{1}{q_t} := \frac{1-t}{q_0}+\frac{t}{q_1}.
\end{gather*}
\end{cor}

\begin{rem}
In the following, we will need to interpolate some inequalities that hold only for functions belonging to a proper subspace of an $L^p$ space, for example only for radial functions, so we notice that if $\mathcal{G}$ is a family of simple functions on $M$ such that the hypotheses of one of the above interpolation theorems hold for $f\in\mathcal{G}$, for instance
\[
	||T_{j+iy} f||_{q_j} \leq M_j(y) ||f||_{p_j} \quad (j=0,1) \quad  \forall f\in \mathcal{G},
\]
then the corresponding interpolation result states that
\[
	||T_t f||_{q_t} \leq M_t ||f||_{p_t} 
\]
for every $f$ in the closure of $\mathcal{G}$ as a subset of $L^{p_t}(\mu)$.
\end{rem}

\begin{lemma} \label{cRub}
If $C(p,s)$ is the optimal constant in Rubin's inequality and $0\leq s \leq s_1 < d/2$, then we have the estimate

\end{lemma}
\[
	C(p_\theta,s) \leq C(p_1, s_1)^{s/s_1} ,
\]
where $\theta=s/s_1$ and $1/p_\theta =(1-\theta)/2+\theta/p_1$.
\begin{proof}
Setting $f = |D|^{-s} \phi$, the Rubin inequality reads
\[
	\left \Vert |\cdot|^{-\beta} \ |D|^{-s} \phi\right \Vert_p \leq C(p,s) ||\phi||_2,
\]
so we are led to consider the family of operators
\[
	\widetilde{T}_{s,\beta} : \phi \mapsto |\cdot|^{-\beta} \ |D|^{-s}\phi .
\]
If $s_1, p_1$ are fixed, to apply the interpolation result we use $\theta=s/s_1$ as a parameter and we set
\[
	\frac{1}{p_\theta} = \frac{1-\theta}{2} + \frac{\theta}{p_1}
\]
and, by the relation between the parameters in Rubin inequality, we have
\[
	\beta_\theta =\theta s_1 +\frac{d}{p_\theta} -\frac{d}{2} = \theta \l \frac{d}{p_1} -\frac{d}{2}+s_1 \r = \theta \beta_1.
\]
Now we have a new family of operators that depend on a single parameter
\[
	T_\theta : \phi \mapsto |\cdot|^{-\beta_\theta} \ |D|^{-s_\theta} \phi = |\cdot|^{-\beta_1 \theta} \ |D|^{-s_1 \theta} \phi
\]
that can be defined on the whole strip $\overline{\Sigma}$. On the left line $\{ x=0\}$ we have
\[
	||T_{iy} \phi||_2 = \lv |D|^{-s_1 iy} \phi \rv_2 = \lv |\cdot|^{-s_1 iy} \h{\phi} \rv_2 = ||\phi||_2
\]
so that
\[
	||T_{iy}||_{L^2 \to L^2 } \leq 1 .
\]
On the left line $\{x=1\}$ 
\begin{align*}
	||T_{1+iy} \phi||_{p_1} &= \lv |\cdot|^{-\beta_1 (1+iy) } |D|^{-s_1 (1+iy)} \phi \rv _{p_1} = \\
	&= \lv |\cdot|^{-\beta_1} |D|^{-s_1} \l |D|^{-s_1 iy} \phi \r \rv_{p_1}
\end{align*}
and applying Rubin's inequality
\[
	||T_{1+iy} \phi||_{p_1} \leq C(p_1,s_1) \lv |D|^{-s_1 iy} \phi \rv_2 = C(p_1,s_1) ||\phi||_2.
\]
Given $\phi,\psi$ simple functions, 
\[
	h(z):=\langle \overline{T_z \phi} ,\psi \rangle = \int_{\R^d} |x|^{-\beta_1 z} \ |D|^{-s_1 z} \phi(x) \ \psi(x) \diff x 
\]
is a continuous function on $\overline{\Sigma}$ and holomorphic in $\Sigma$; indeed consider the two variable function
\begin{align*}
	H(z_1,z_2) &:= \int_{\R^d} |x|^{-\beta_1 z_1} \ |D|^{-s_1 z_2} \phi(x) \ \psi(x) \diff x = \\
		&= \int_{\R^d} \mathcal{F}\left[ |x|^{-\beta_1 z_1} \psi(x)\right](\xi) \ |\xi|^{-s_1z_2} \h{\phi}(\xi)  \diff \xi.
\end{align*}
The first expression for $H$ shows that, fixing $z_2$, $\frac{\partial H(z_1,z_2)}{\partial \bar z_1}=0$, while the second form gives, $z_1$ being fixed, $\frac{\partial H(z_1,z_2)}{\partial \bar z_2}=0$. Hence $H$ is holomorphic in $\Sigma\times \Sigma$ (Osgood's Lemma) and 
\[
	z\mapsto H(z,z)=h(z)
\]
is a composition of holomorphic functions, so it is holomorphic itself. Moreover, the function $h$ is bounded on the strip, since
\begin{align*}
	|h(z)| &= |\langle |D|^{-s_1 \bar{z}} \overline{\phi} , |\cdot|^{-\beta_1 \bar{z}} \psi \rangle | 
	= | \langle |\cdot|^{-s_1 \bar{z}} \h{\overline{\phi}} , \mathcal{F} (|\cdot|^{-\beta_1 \bar{z}} \psi)) \rangle | \leq \\
	&\leq| \langle |\cdot|^{-s_1 x} \h{\overline{\phi}} , \mathcal{F} (|\cdot|^{-\beta_1 \bar{z}} \psi)) \rangle | 
	= |\langle |D|^{-s_1 x} \overline{\phi} , |\cdot|^{-\beta_1\bar{z}} \psi \rangle | \leq \\
	&\leq |\langle |D|^{-s_1 x} \overline{\phi} , |\cdot|^{-\beta_1 x} \psi \rangle | = |\langle \overline{T_x \phi} ,\psi \rangle=|h(x)|,
\end{align*}
where we have used repeatedly Plancherel identity. This proves that $|h|$ is bounded on the strip if and only if it is bounded in $[0,1]$; by continuity, the function $|h|$ is indeed bounded there. Thus we apply Stein interpolation to get
\[
	||T_\theta \phi||_{p_\theta} \leq \l C(p_1,s_1)\r^\theta ||\phi||_2 ,
\] 
which means
\[
	\lv |\cdot|^{-\beta_1 s/s_1} |D|^{-s} \phi \rv_{p_\theta} \leq \l C(p_1,s_1)\r^\theta ||\phi||_2 ,
\]
for every $\phi \in L_{rad}^2(\R^d)$ and substituting back $\phi = |D|^s f$
\[
	\lv |\cdot|^{-\beta_1 s/s_1} f \rv_{p_\theta} \leq  \l C(p_1,s_1)\r^\theta ||f||_{\dot{H}^s}
\]
for every $f\in |D|^{-s}(L_{rad}^2(\R^d))= \dot{H}_{rad}^s(\R^d)$.
\end{proof}

\begin{proof}[Proof of Main Theorem]
We start by differentiating all the relevant quantities: for the weighted $L^p$ norm
\begin{gather*}
	\partial_p \lv  |\cdot|^{-\beta} f \rv_p = \partial_p \l  \lv |\cdot|^{-\beta} f \rv_p^p \r^{1/p} = \\
	=\lv  |\cdot|^{-\beta} f \rv_p  \partial_p\l \frac{1}{p} \log \l \lv |\cdot|^{-\beta} f \rv_p^p \r \r = \\
	= \lv  |\cdot|^{-\beta} f \rv_p \l -\frac{1}{p^2} \log \l \lv  |\cdot|^{-\beta} f \rv_p^p \r 
	+ \frac{1}{p} \frac{\int_{\R^d} \partial_p \l |f|^p |x|^{-\beta p} \r \diff x}{\lv  |\cdot|^{-\beta} f \rv_p^p} \r 
\end{gather*}
and
\begin{gather*}
	\partial_{\beta} \lv  |\cdot|^{-\beta} f \rv_p = \partial_\beta \l  \lv |\cdot|^{-\beta} f \rv_p^p \r^{1/p} = \\
	= \frac{1}{p} \l  \lv |\cdot|^{-\beta} f \rv_p^p \r^{1/p-1} \partial_\beta \l \lv |\cdot|^{-\beta} f \rv_p^p \r = \\
	=  \frac{1}{p} \l  \lv |\cdot|^{-\beta} f \rv_p^p \r^{1/p-1} \int_{\R^d}  |f|^p |x|^{-\beta p} \l -p \log |x| \r \diff x .
\end{gather*}
Computing the above quantities for $p=2, \beta=0$ 
\begin{gather*}
	\partial_p \lv  |\cdot|^{-\beta} f \rv_p \Big|_{\substack{ {p=2} \\ {\beta=0} }} 
	= ||f||_2 \l -\frac{1}{4} \log ||f||_2^2 + \frac{1}{2} ||f||_2^{-2} \int_{\R^d} |f|^2 \log|f|  \diff x  \r \\
	= \frac{1}{2} \int_{\R^d} |f|^2 \log |f| \diff x 
\end{gather*}
and similarly
\[
	\partial_\beta \lv  |\cdot|^{-\beta} f \rv_p \Big|_{\substack{ {p=2} \\ {\beta=0} }}  
	= - ||f||_2^{-1} \! \int_{\R^d} |f|^2 \log |x| \diff x =  - \int_{\R^d} |f|^2 \log |x| \diff x .
\]
This means that
\begin{gather*}
	\frac{\diff}{\diff s} \lv  |\cdot|^{-\beta} f \rv_p |_{s=0}
	= \frac{\diff p}{\diff s} |_{s=0} \ \partial_p \lv  |\cdot|^{-\beta} f \rv_p \Big|_{\substack{ {p=2} \\ {\beta=0} }} + \\
	+\frac{\diff \beta}{\diff s} |_{s=0} \ \partial_\beta \lv  |\cdot|^{-\beta} f \rv_p \Big|_{\substack{ {p=2} \\ {\beta=0} }}  = \\
	= \frac{1}{2} \frac{\diff p}{\diff s} |_{s=0} \int_{\R^d} |f|^2 \log |f| \diff x
	- \frac{\diff \beta}{\diff s} |_{s=0} \int_{\R^d} |f|^2 \log |x| \diff x  .
\end{gather*}
Now we choose $s_1, p_1, \beta_1$ in the permitted range of values and, referring to \cref{cRub}, define
\[
	p(s) := p_{s/s_1} \quad \text{ and } \quad \beta(s) := \beta_{s/s_1} .
\]
With this choice 
\begin{gather*}
	\frac{\diff p}{\diff s} |_{s=0} = \frac{4}{s_1} \l \frac{1}{2} - \frac{1}{p_1} \r = \frac{4}{d} \l 1-\frac{\beta_1}{s_1} \r \\
	\frac{\diff \beta}{\diff s} |_{s=0}  = 1-\frac{d}{s_1} \l \frac{1}{2}-\frac{1}{p_1} \r = \frac{\beta_1}{s_1} .
\end{gather*}
Since we are considering the Rubin inequality with a constant given by \cref{cRub}, we need to compute
\[
	\partial_s \l C(p_1,s_1) \r^{s/s_1} |_{s=0} = \frac{1}{s_1} \log C(p_1,s_1)
\]
and recalling that
\[
	\partial_s ||f||_{\dot{H}^s} |_{s=0} =  \int_{\R^d} |\h{f}|^2 \log |\xi| \diff \xi 
\]
we get the inequality
\begin{equation}
\begin{split}
	\frac{2}{d}\l 1-\frac{\beta_1}{s_1}\r \int_{\R^d} |f|^2 \log |f| \diff x
	- \frac{\beta_1}{s_1} \int_{\R^d} |f|^2 \log |x| \diff x   \leq \\
	\leq \frac{1}{s_1} \log C(p_1,s_1) + \int_{\R^d} |\h{f}|^2 \log |\xi| \diff \xi  .
\end{split}
\end{equation}
In particular we can choose the values $p_1=2/(1-2s_1)$ and $\beta_1=s_1(1-d)$ to get
\[
	2 \int_{\R^d} |f|^2 \log |f| \diff x + (d-1) \int_{\R^d} |f|^2 \log |x| \diff x  \leq  \int_{\R^d} |\h{f}|^2 \log |\xi| \diff \xi  +c
\]
for some $c\in\R$.
\end{proof}

\bibliographystyle{plain}
\bibliography{Bibliografia}
\end{document}